\documentclass{article}
\usepackage{geometry}    
\usepackage{amsfonts}    
\usepackage{amssymb}     
\usepackage[all]{xy}     
\usepackage{amsmath}     
\usepackage[dvips]{graphicx}       
\usepackage{color}
\usepackage{pst-tree}    
\usepackage{pstricks}

\newcommand{\go}[1]{\mathfrak {#1}}
\newcommand{\rrbr}{]\!]}
\newcommand{\llbr}{[\![}

\newtheorem{The}{Theorem}

\newtheorem{Pro}[The]{Proposition}

\newtheorem{Rk}[The]{Remark}

\newtheorem{Conj}[The]{Conjecture}

\begin{document}

\begin{center}
\Large{Branching Law for the Finite Subgroups of $\mathbf{SL}_4\mathbb{C}$}
\end{center}

\vspace{.4cm}
\begin{center}
Fr\'ed\'eric BUTIN\footnote{Universit\'e de Lyon, Universit\'e Lyon 1, CNRS, UMR5208, Institut Camille Jordan, 43 blvd du 11 novembre 1918, F-69622 Villeurbanne-Cedex, France, email: butin@math.univ-lyon1.fr}\vspace{.5cm}\\
\end{center}

\begin{small}
\noindent \textbf{\textsc{Abstract}}\\
In the framework of McKay correspondence we determine, for every finite subgroup $\Gamma$ of $\mathbf{SL}_4\mathbb{C}$, how the finite dimensional irreducible representations of $\mathbf{SL}_4\mathbb{C}$ decompose under the action of $\Gamma$.\\
Let $\go{h}$ be a Cartan subalgebra of $\go{sl}_4\mathbb{C}$ and let $\varpi_1,\,\varpi_2,\,\varpi_3$ be the corresponding fundamental weights. For $(p,q,r)\in \mathbb{N}^3$, the restriction $\pi_{p,q,r}|_\Gamma$ of the irreducible representation $\pi_{p,q,r}$ of highest weight $p\varpi_1+q\varpi_2+r\varpi_3$ of $\mathbf{SL}_4\mathbb{C}$ decomposes as ${\pi_{p,q,r}}|_{\Gamma}=\bigoplus_{i=0}^l m_i(p,q,r)\gamma_i.$ We determine the multiplicities $m_i(p,q,r)$ and prove that the series $P_\Gamma(t,u,w)_i=\sum_{p=0}^\infty\sum_{q=0}^\infty\sum_{r=0}^\infty m_i(p,q,r)t^pu^qw^r$ are rational functions.\\
This generalizes results from Kostant for $\mathbf{SL}_2\mathbb{C}$ and our preceding works about $\mathbf{SL}_3\mathbb{C}$.\\
\end{small}

\begin{small}

\end{small}

\section{\textsf{Introduction and results}}\label{firstsection}

\noindent $\bullet$ Let $\Gamma$ be a finite subgroup of
$\mathbf{SL}_4\mathbb{C}$ and
$\{\gamma_0,\dots,\,\gamma_l\}$ the set of
equivalence classes of irreducible finite dimensional complex
representations of $\Gamma$, where $\gamma_0$ is the trivial
representation. The character
associated to $\gamma_j$ is denoted by $\chi_j$.\\
Consider $\gamma : \Gamma\rightarrow
\mathbf{SL}_4\mathbb{C}$ the natural $4-$dimensional
representation, and $\gamma^*$ its contragredient representation. The character of $\gamma$ is denoted by $\chi$.
By complete reducibility  we get the decompositions
$$\forall\ j\in\llbr 0,\,l\rrbr,\,
\ \gamma_j\otimes\gamma=\bigoplus_{i=0}^l a^{(1)}_{ij}\gamma_i\ \
,\,
\ \gamma_j\otimes(\gamma\wedge\gamma)=\bigoplus_{i=0}^l a^{(2)}_{ij}\gamma_i\ \
\textrm{and}\ \ \gamma_j\otimes\gamma^*=\bigoplus_{i=0}^l
a^{(3)}_{ij}\gamma_i.$$
This defines the three following square matrices of $\mathbf{M}_{l+1}\mathbb{N}$:
$$A^{(1)}:=\left(a^{(1)}_{ij}\right)_{(i,j)\in\llbr 0,\, l\rrbr^2},\ A^{(2)}:=\left(a^{(2)}_{ij}\right)_{(i,j)\in\llbr 0,\, l\rrbr^2}\ \textrm{and}\ A^{(3)}:=\left(a^{(3)}_{ij}\right)_{(i,j)\in\llbr 0,\,l\rrbr^2}.$$

\noindent $\bullet$ Let $\go{h}$ be a Cartan subalgebra of $\go{sl}_4\mathbb{C}$ and let $\varpi_1,\ \varpi_2,\ \varpi_3$ be the corresponding fundamental weights, and $V(p\varpi_1+q\varpi_2+r\varpi_3)$ the simple
$\go{sl}_4\mathbb{C}-$module of highest weight $p\varpi_1+q\varpi_2+r\varpi_3$ with $(p,\,q,\,r)\in \mathbb{N}^3$. Then we get an irreducible representation
$\pi_{p,q,r}:\mathbf{SL}_4\mathbb{C}\rightarrow
\mathbf{GL}(V(p\varpi_1+q\varpi_2+r\varpi_3))$.
 The restriction of
$\pi_{p,q,r}$ to the subgroup $\Gamma$ is a representation
of~$\Gamma$, and by complete reducibility, we get the
decomposition $${\pi_{p,q,r}}|_{\Gamma}=\bigoplus_{i=0}^l
m_i(p,q,r)\gamma_i,$$ where the $m_i(p,q,r)$'s are non negative integers.
Let $\mathcal{E}:=(e_0,\dots,\,e_l)$ be the canonical basis of
$\mathbb{C}^{l+1}$, and $$\displaystyle{v_{p,q,r}:=\sum_{i=0}^l
m_i(p,q,r)e_i\in \mathbb{C}^{l+1}}.$$\\
We have in particular $v_{0,0,0}=e_0$ as $\gamma_0$ is the trivial representation. Let us consider the vector
$$\displaystyle{P_\Gamma(t,\,u,\,w):=\sum_{p=0}^\infty\sum_{q=0}^\infty\sum_{r=0}^\infty
v_{p,q,r}t^pu^qw^r\in\left(\mathbb{C}\llbr t,\,u,\,w\rrbr\right)^{l+1}},$$
and denote by $\displaystyle{P_\Gamma(t,\,u,\,w)_j}$ its $j-$th coordinate in the basis $\mathcal{E}$, which is an element of $\mathbb{C}\llbr t,\,u,\,w\rrbr$.
Note that $P_\Gamma(t,\,u,\,w)$ can also be seen as a formal power
series with coefficients in $\mathbb{C}^{l+1}$. The aim of this article is to prove the following theorem.

\begin{The}\label{mainth}$\\$
The coefficients of $P_\Gamma(t,\,u,\,w)$ are rational fractions in $t,\,u,\,w$, i.e. the formal power series $P_\Gamma(t,\,u,\,w)_i$ are rational functions $$P_\Gamma(t,\,u,\,w)_i=\frac{N_\Gamma(t,\,u,\,w)_i}{D_\Gamma(t,\,u,\,w)},\ i\in\llbr 0,\,l\rrbr,$$
where the $N_\Gamma(t,\,u,\,w)_i$'s and $D_\Gamma(t,\,u,\,w)$ are elements of $\mathbb{Q}[t,\,u,\,w]$.\\
\end{The}

\noindent $\bullet$ The proof of this theorem uses a key-relation satisfied by $P_\Gamma(t,\,u,\,w)$ as well as a so-called inversion formula. Two essential ingredients are the decomposition of the tensor product of $\pi_{p,q,r}$ with the natural representation of $\mathbf{SL}_4\mathbb{C}$ and the simultaneous diagonalizability of certain matrices. The effective calculation of $P_\Gamma(t,\,u,\,w)$ then reduces to matrix multiplication.\\

\noindent In \cite{BP09} we applied a similar method for $\mathbf{SL}_2\mathbb{C}$ --- recovering thereby in a quite easy way the results obtained by Kostant in \cite{Kos85}, \cite{Kos06}, and by Gonzalez-Sprinberg and Verdier in \cite{GSV83} --- and for $\mathbf{SL}_3\mathbb{C}$ in order to get explicit computations of the series for every finite subgroup of $\mathbf{SL}_3\mathbb{C}$.\\

\noindent The general framework of that study is the construction of a minimal resolution of singularities of the orbifold $\mathbb{C}^n/\Gamma$. It is related to the McKay correspondence (see \cite{BKR01}, \cite{GSV83} and \cite{GNS04}). For example, Gonzalez-Spriberg and Verdier use in \cite{GSV83} a Poincar\'e series to construct explicitly minimal resolutions for singularities of $V=\mathbb{C}^2/\Gamma$ when $\Gamma$ is a finite subgroup of $\mathbf{SL}_2\mathbb{C}$. To go further in this approach, our results for $\mathbf{SL}_4\mathbb{C}$ could be used to construct an explicit synthetic minimal resolution of singularities for orbifolds of the form $\mathbb{C}^4/\Gamma$ where $\Gamma$ is a finite subgroup of $\mathbf{SL}_4\mathbb{C}$.

\section{\textsf{Properties of the matrices $A^{(1)},\,A^{(2)},\,A^{(3)}$}}

\noindent In order to compute the series $P_\Gamma(t,\,u,\,w)$, we first establish here some properties of the matrices $A^{(1)},\,A^{(2)},\,A^{(3)}$. The first proposition essentially follows from the uniqueness of the decomposition of a representation as sum of irreducible representations.

\begin{Pro}\label{propmatrices}$\\$
$\bullet$ $A^{(3)}=\,^tA^{(1)}$.\\
$\bullet$ $A^{(2)}$ is a symmetric matrix.\\
$\bullet$ $A^{(1)},\ A^{(2)}$ and $A^{(3)}$ commute. In particular, $A^{(1)}$ is a normal matrix.
\end{Pro}

\underline{Proof:}\\
Since $a^{(1)}_{ij}=(\chi_i\,|\,\chi_{\gamma\otimes\gamma_j})=\frac{1}{|\Gamma|}\sum_{g\in\Gamma}\overline{\chi_i(g)}\chi(g)\chi_j(g)$, we have $\gamma\otimes\gamma_j=\bigoplus_{i=0}^l a^{(1)}_{ij}\gamma_i$. In the same way, $(\gamma\wedge\gamma)\otimes\gamma_j=\bigoplus_{i=0}^l a^{(2)}_{ij}\gamma_i$ and $\gamma^*\otimes\gamma_j=\bigoplus_{i=0}^l a^{(3)}_{ij}\gamma_i$.\\
Then
$$\begin{array}{rcl}
a^{(3)}_{ij} & = & (\chi_i\,|\,\chi_{\gamma_j\otimes\gamma^*})=\frac{1}{|\Gamma|}\sum_{g\in\Gamma}\overline{\chi_i(g)}\chi_j(g)\chi_{\gamma^*}(g)
=\frac{1}{|\Gamma|}\sum_{g\in\Gamma}\overline{\chi_i(g)}\chi_j(g)\chi(g^{-1})\\
 & = & \frac{1}{|\Gamma|}\sum_{g\in\Gamma}\overline{\chi_i(g^{-1})}\chi_j(g^{-1})\chi(g)
 =\frac{1}{|\Gamma|}\sum_{g\in\Gamma}\chi_i(g)\overline{\chi_j(g)}\chi(g)=a^{(1)}_{ji},
      \end{array}$$
hence $A^{(3)}=\,^tA^{(1)}$.\\
We also have $(\gamma_j\otimes \gamma)\otimes \gamma^*=\left(\bigoplus_{i=0}^l a^{(1)}_{ij}\gamma_i\right)\otimes\gamma^*
=\bigoplus_{i=0}^l a^{(1)}_{ij}\left(\bigoplus_{k=0}^l a^{(3)}_{ki}\gamma_k\right)
=\bigoplus_{k=0}^l\left(\sum_{i=0}^l a^{(3)}_{ki}a^{(1)}_{ij}\right)\gamma_k$
and $\gamma\otimes(\gamma_j\otimes \gamma^*)=\gamma\otimes\left(\bigoplus_{i=0}^l a^{(3)}_{ij}\gamma_i\right)
=\bigoplus_{i=0}^l a^{(3)}_{ij}\left(\bigoplus_{k=0}^l a^{(1)}_{ki}\gamma_k\right)
=\bigoplus_{k=0}^l\left(\sum_{i=0}^l a^{(1)}_{ki}a^{(3)}_{ij}\right)\gamma_k$, hence $A^{(3)}A^{(1)}=A^{(1)}A^{(3)}$. The proofs of the other statements are the same. $\blacksquare$\\

\noindent Since $A^{(1)},\ A^{(2)},\ A^{(3)}$ are normal, we know that they are
diagonalizable with eigenvectors forming an orthogonal basis. Now we will
diagonalize these matrices by using the character table of
the group $\Gamma$.
Let us denote by $\{C_0,\dots,\,C_l\}$ the set of conjugacy
classes of $\Gamma$, and for any $j\in\llbr 0,\,l\rrbr$, let $g_j$
be an element of $C_j$. So the character table of $\Gamma$ is the
matrix $T_\Gamma\in \mathbf{M}_{l+1}\mathbb{C}$ defined by~$(T_\Gamma)_{i,j}:=\chi_i(g_j)$.\\

\begin{Pro}\label{vp}$\\$
$\bullet$ For $k\in\llbr 0,\,l\rrbr$, set
$w_k:=(\chi_0(g_k),\dots,\,\chi_l(g_k))\in\mathbb{C}^{l+1}$. Then
$w_k$ is an eigenvector of $A^{(3)}$ associated to the eigenvalue
$\chi(g_k)$. Similarly, $w_k$ is an eigenvector of $A^{(1)}$ associated to the
eigenvalue $\overline{\chi(g_k)}$.\\
$\bullet$ For $k\in\llbr 0,\,l\rrbr$, $w_k$ is an eigenvector of $A^{(2)}$ associated to the eigenvalue
$\frac{1}{2}\left(\chi(g_k)^2+\chi(g_k^2)\right)$.
\end{Pro}

\underline{Proof:}\\
From the relation $\gamma_i \otimes\gamma=\sum_{j=0}^l a^{(1)}_{ji}\gamma_j$, we get $\chi_i\chi=\chi_{\gamma_i \otimes\gamma}=\sum_{j=0}^l a^{(1)}_{ji}\chi_j$. By evaluating this on $g_k$, we obtain $\chi_i(g_k)\chi(g_k)=\sum_{j=0}^l a^{(1)}_{ji}\chi_j(g_k)=\sum_{j=0}^l a^{(3)}_{ij}\chi_j(g_k)$ according to Proposition \ref{propmatrices}. So $w_k$ is an eigenvector of $A^{(3)}$ associated to the eigenvalue $\chi(g_k)$. The method is similar for the other results. $\blacksquare$\\

\noindent As the $w_j$'s are the column of $T_\Gamma$, which are always orthogonal, the matrix $T_\Gamma$ is invertible and the family $\mathcal{W}:=(w_0,\dots,\,w_l)$ is a common basis of eigenvectors of $A^{(1)}$, $A^{(2)}$ and $A^{(3)}$. Then $\Lambda^{(1)}:=T_\Gamma^{-1}\,A^{(1)}\,T_\Gamma$, $\Lambda^{(2)}:=T_\Gamma^{-1}\,A^{(2)}\,T_\Gamma$ and $\Lambda^{(3)}:=T_\Gamma^{-1}\,A^{(3)}\,T_\Gamma$
are diagonal matrices, with $\Lambda^{(1)}_{jj}=\overline{\chi(g_j)}$, $\Lambda^{(2)}_{jj}=\frac{1}{2}(\chi(g_j)^2-\chi(g_j^2))$ and $\Lambda^{(3)}_{jj}=\chi(g_j)$.\\

\noindent Now, we make use of the Clebsch-Gordan formula
\begin{equation}\label{CG}
    \begin{array}{rcl}
      \pi_{1,0,0}\otimes\pi_{p,q,r} & = & \pi_{p+1,q,r}\oplus\pi_{p,q,r-1}\oplus\pi_{p-1,q+1,r}\oplus\pi_{p,q-1,r+1}, \\
      \pi_{0,1,0}\otimes\pi_{p,q,r} & = & \pi_{p,q+1,r}\oplus\pi_{p,q-1,r}\oplus\pi_{p+1,q-1,r+1}\oplus\pi_{p-1,q+1,r-1}
      \oplus\pi_{p-1,q,r+1}\oplus\pi_{p+1,q,r-1}, \\
      \pi_{0,0,1}\otimes\pi_{p,q,r} & = & \pi_{p,q,r+1}\oplus\pi_{p-1,q,r}\oplus\pi_{p,q+1,r-1}\oplus\pi_{p+1,q-1,r}. \\
    \end{array}
\end{equation}

\begin{Pro}$\\$
The vectors $v_{m,n}$ satisfy the following recurrence relations
$$\begin{array}{rcl}
  A^{(1)}v_{p,q,r} & = & v_{p+1,q,r}+v_{p,q,r-1}+v_{p-1,q+1,r}+v_{p,q-1,r+1}, \\
  A^{(2)}v_{p,q,r} & = & v_{p,q+1,r}+v_{p,q-1,r}+v_{p+1,q-1,r+1}+v_{p-1,q+1,r-1}
      +v_{p-1,q,r+1}+v_{p+1,q,r-1}, \\
  A^{(3)}v_{p,q,r} & = & v_{p,q,r+1}+v_{p-1,q,r}+v_{p,q+1,r-1}+v_{p+1,q-1,r}.
\end{array}$$
\end{Pro}

\underline{Proof:}\\
The definition of $v_{p,q,r}$ reads
$v_{p,q,r}=\sum_{i=0}^l m_i(p,\,q,\,r)e_i$, thus
$A^{(1)}v_{p,q,r}=\sum_{i=0}^l\left(\sum_{j=0}^l
m_j(p,\,q,\,r)a^{(1)}_{ij}\right)e_i$.\\
Now
$$(\pi_{1,0,0}\otimes\pi_{p,q,r})|_\Gamma=\pi_{p,q,r}|_\Gamma\otimes\gamma=\sum_{j=0}^l
m_j(p,\,q,\,r)\gamma_j\otimes\gamma=\sum_{i=0}^l \left(\sum_{j=0}^l
m_j(p,\,q,\,r)a^{(1)}_{ij}\right)\gamma_i,$$and
$$\begin{array}{l}
\displaystyle{\pi_{p+1,q,r}|_\Gamma+\pi_{p,q,r-1}|_\Gamma+\pi_{p-1,q+1,r}|_\Gamma+\pi_{p,q-1,r+1}|_\Gamma} \\
\ \ \ \ \ \ \ \ \ \ \ \ \ \ \ \ \ = \displaystyle{\sum_{i=0}^l
\left(m_i(p+1,\,q,\,r)+m_i(p,\,q,\,r-1)+m_i(p-1,\,q+1,\,r)+m_i(p,\,q-1,\,r+1)\right)\gamma_i}.
  \end{array}
$$
By uniqueness, $$\sum_{j=0}^l
m_j(p,\,q,\,r)a^{(1)}_{ij}=m_i(p+1,\,q,\,r)+m_i(p,\,q,\,r-1)+m_i(p-1,\,q+1,\,r)+m_i(p,\,q-1,\,r+1).\ \blacksquare$$

\section{\textsf{The series $P_\Gamma(t,\,u,\,w)$ is a rational function}}

\noindent This section is mainly devoted to the proof of Theorem \ref{mainth}.

\subsection{\textsf{A key-relation satisfied by the series $P_\Gamma(t,\,u,\,w)$}}

\begin{Pro}\label{form}$\\$
Set
$$\begin{array}{rcl}
J(t,\,u,\,w) & := & (1-u^2)((1+ut^2)(1+uw^2)-tw(1+u^2))I_n+twu(1-u^2)A^{(2)}\\
 & & -tu(1+uw^2)(A^{(3)}-uA^{(1)})-wu(1+ut^2)(A^{(1)}-uA^{(3)}).
      \end{array}$$
Then the series $P_\Gamma(t,\,u,\,w)$ satisfies the following relation
$$\begin{array}{rcl}
    J(t,\,u,\,w)\,v_{0,0,0} & = & \left(1-tA^{(1)}+t^2A^{(2)}-t^3A^{(3)}+t^4\right)\left(1-wA^{(3)}+w^2A^{(2)}-w^3A^{(1)}+w^4\right) \\
 & & \Big((1+u^2)(1-u^2)^2-u(1-u^2)^2A^{(2)}+u^2(A^{(1)}-uA^{(3)})(A^{(3)}-uA^{(1)})\Big)P_\Gamma(t,\,u,\,w).
  \end{array}$$
\end{Pro}

\underline{Proof:}\\
$\bullet$ Set $x:=P_\Gamma(t,\,u,\,w)$. Set also $v_{p,q,-1}:=0$, $v_{p,-1,r}:=0$ and
$v_{-1,q,r}:=0$ for $(p,\,q,\,r)\in\mathbb{N}^3$, such that, according to
the Clebsch-Gordan formula, the formulae of the preceding
corollary are still true for~$(p,\,q,\,r)\in\mathbb{N}^3$. So we have (by denoting $\sum_{p=0}^\infty\sum_{q=0}^\infty\sum_{r=0}^\infty$ by $\sum_{pqr}$)
$$\begin{array}{rl}
&(1-wA^{(3)}+w^2A^{(2)}-w^3A^{(1)}+w^4)x \\ \\
= & \displaystyle{\sum_{pqr}v_{p,q,r}t^pu^qw^r-\sum_{pqr}(v_{p,q,r+1}+v_{p-1,q,r}+v_{p,q+1,r-1}+v_{p+1,q-1,r})t^pu^qw^{r+1}} \\
       & \displaystyle{+\sum_{pqr}(v_{p,q+1,r}+v_{p,q-1,r}+v_{p+1,q-1,r+1}+v_{p-1,q+1,r-1}+v_{p-1,q,r+1}+v_{p+1,q,r-1})t^pu^qw^{r+2}} \\
        & \displaystyle{-\sum_{pqr}(v_{p+1,q,r}+v_{p,q,r-1}+v_{p-1,q+1,r}+v_{p,q-1,r+1})t^pu^qw^{r+3}+\sum_{pqr}v_{p,q,r}t^pu^qw^{r+4}},
  \end{array}
$$
hence
\begin{equation}\label{rel1}
    (1-wA^{(3)}+w^2A^{(2)}-w^3A^{(1)}+w^4)x=(1-tw+uw^2-t^{-1}uw)\sum_{p=0}^\infty\sum_{q=0}^\infty v_{p,q,0}t^pu^q+t^{-1}uw\sum_{q=0}^\infty v_{0,q,0}u^q.
\end{equation}
$\bullet$ In the same way (by denoting $\sum_{p=0}^\infty\sum_{q=0}^\infty$ by $\sum_{pq}$)
$$\begin{array}{rl}
&\displaystyle{(1-tA^{(1)}+t^2A^{(2)}-t^3A^{(3)}+t^4)\sum_{p=0}^\infty\sum_{q=0}^\infty v_{p,q,0}t^pu^q} \\ \\
= & \displaystyle{\sum_{pq}v_{p,q,0}t^pu^q-\sum_{pq}(v_{p+1,q,0}+v_{p-1,q+1,0}+v_{p,q-1,1})t^{p+1}u^q} \\
       & \displaystyle{+\sum_{pq}(v_{p,q+1,0}+v_{p,q-1,0}+v_{p+1,q-1,1}+v_{p-1,q,1})t^{p+2}u^q} \\
        & \displaystyle{-\sum_{pq}(v_{p,q,1}+v_{p-1,q,0}+v_{p+1,q-1,0})t^{p+3}u^q+\sum_{pq}v_{p,q,0}t^{p+4}u^q},
  \end{array}
$$
hence
\begin{equation}\label{rel2}
    (1-tA^{(1)}+t^2A^{(2)}-t^3A^{(3)}+t^4)\sum_{p=0}^\infty\sum_{q=0}^\infty v_{p,q,0}t^pu^q=(1+t^2u)\sum_{q=0}^\infty v_{0,q,0}u^q-tu\sum_{q=0}^\infty v_{0,q,1}u^q.
\end{equation}
Moreover, we have
$$\begin{array}{rl}
& \displaystyle{(1-tA^{(1)}+t^2A^{(2)}-t^3A^{(3)}+t^4)\sum_{q=0}^\infty v_{0,q,0}u^q} \\
= & \displaystyle{\sum_{q=0}^\infty v_{0,q,0}u^q-\sum_{q}^\infty(v_{1,q,0}+v_{0,q-1,1})tu^q}\displaystyle{+\sum_{q}^\infty(v_{0,q+1,0}+v_{0,q-1,0}+v_{1,q-1,1})t^2u^q} \\ & \displaystyle{-\sum_{q}^\infty(v_{0,q,1}+v_{1,q-1,0})t^{3}u^q+\sum_{q}^\infty v_{0,q,0}t^{4}u^q},
  \end{array}
$$
hence
\begin{equation}\label{rel3}
\begin{array}{rl}
 & \displaystyle{(1-tA^{(1)}+t^2A^{(2)}-t^3A^{(3)}+t^4)\sum_{p=0}^\infty v_{0,q,0}u^q}\\
 = & \displaystyle{(1+t^4+t^2u^{-1}+t^2u)\sum_{q=0}^\infty v_{0,q,0}u^q-t^2u^{-1}v_{0,0,0}-(t+t^3u)\sum_{q=0}^\infty v_{1,q,0}u^q}\\
  & \displaystyle{-(tu+t^3)\sum_{q=0}^\infty v_{0,q,1}u^q+t^2u\sum_{q=0}^\infty v_{1,q,1}u^q}.
\end{array}
\end{equation}
By combining Equations (\ref{rel1}), (\ref{rel2}) and (\ref{rel3}), we get
$$\begin{array}{rl}
 & (1-tA^{(1)}+t^2A^{(2)}-t^3A^{(3)}+t^4)(1-wA^{(3)}+w^2A^{(2)}-w^3A^{(3)}+w^4)x \\ \\
= & \displaystyle{(1-tA^{(1)}+t^2A^{(2)}-t^3A^{(3)}+t^4)\left((1-tw+uw^2-t^{-1}uw)\sum_{pq}v_{p,q,0}t^pu^q+t^{-1}uw\sum_{q=0}^\infty v_{0,q,0}u^q\right)}\\ \\
= &  \displaystyle{(1-tw+uw^2-t^{-1}uw)\left((1+t^2u)\sum_{q=0}^\infty v_{0,q,0}u^q-tu\sum_{q=0}^\infty v_{0,q,1}u^q\right)}\\
 & \displaystyle{+(1+t^4+t^2u^{-1}+t^2u)t^{-1}uw\sum_{q=0}^\infty v_{0,q,0}u^q-twv_{0,0,0}-(1+t^2u)uw\sum_{q=0}^\infty v_{1,q,0}u^q}\\
 & \displaystyle{-(u+t^2)uw\sum_{q=0}^\infty v_{0,q,1}u^q+tu^2w\sum_{q=0}^\infty v_{1,q,1}u^q},
  \end{array}
$$
hence
\begin{equation}\label{rel4}
\begin{array}{rl}
 & \displaystyle{(1-tA^{(1)}+t^2A^{(2)}-t^3A^{(3)}+t^4)(1-wA^{(3)}+w^2A^{(2)}-w^3A^{(1)}+w^4)x}\\ \\
 = & \displaystyle{(1+ut^2)(1+uw^2)\sum_{q=0}^\infty v_{0,q,0}u^q-tu(1+uw^2)\sum_{q=0}^\infty v_{0,q,1}u^q}\\
  & \displaystyle{-wu(1+ut^2)\sum_{q=0}^\infty v_{1,q,0}u^q-twv_{0,0,0}+tu^2w\sum_{q=0}^\infty v_{1,q,1}u^q}.
\end{array}
\end{equation}
Besides, we have the two following equations
\begin{equation}\label{R1}
    A^{(1)}\sum_{q=0}^\infty v_{0,q,0}u^q=\sum_{q=0}^\infty v_{1,q,0}u^q+u\sum_{q=0}^\infty v_{0,q,1}u^q,
\end{equation}
and
\begin{equation}\label{R3}
    A^{(3)}\sum_{q=0}^\infty v_{0,q,0}u^q=\sum_{q=0}^\infty v_{0,q,1}u^q+u\sum_{q=0}^\infty v_{1,q,0}u^q.
\end{equation}
From these two equations, we deduce
\begin{equation}\label{R7}
    \sum_{q=0}^\infty v_{0,q,1}u^q=(1-u^2)^{-1}(A^{(3)}-uA^{(1)})\sum_{q=0}^\infty v_{0,q,0}u^q.
\end{equation}
Now, we have
\begin{equation}\label{R4}
    A^{(1)}\sum_{q=0}^\infty v_{0,q,1}u^q=\sum_{q=0}^\infty v_{1,q,1}u^q+\sum_{q=0}^\infty v_{0,q,0}u^q+u\sum_{q=0}^\infty v_{0,q,2}u^q,
\end{equation}
and
\begin{equation}\label{R6}
    A^{(3)}\sum_{q=0}^\infty v_{0,q,1}u^q=\sum_{q=0}^\infty v_{0,q,2}u^q+u^{-1}\sum_{q=0}^\infty v_{0,q,0}u^q+u\sum_{q=0}^\infty v_{1,q,1}u^q-u^{-1}v_{0,0,0},
\end{equation}
therefore $$\sum_{q=0}^\infty v_{1,q,1}u^q=(1-u^2)^{-1}(A^{(1)}-uA^{(3)})\sum_{q=0}^\infty v_{0,q,1}u^q-(1-u^2)^{-1}v_{0,0,0}.$$
So, according to Equation (\ref{R7}), we deduce
\begin{equation}\label{R8}
\sum_{q=0}^\infty v_{1,q,1}u^q=(1-u^2)^{-2}(A^{(1)}-uA^{(3)})(A^{(3)}-uA^{(1)})\sum_{q=0}^\infty v_{0,q,0}u^q-(1-u^2)^{-1}v_{0,0,0}.
\end{equation}
By using Equation (\ref{R8}), we may write Equation (\ref{rel4}) as
\begin{equation}\label{rel4b}
\begin{array}{rl}
 & \displaystyle{(1-tA^{(1)}+t^2A^{(2)}-t^3A^{(3)}+t^4)(1-wA^{(3)}+w^2A^{(2)}-w^3A^{(1)}+w^4)x}\\ \\
 = & \displaystyle{\left((1+ut^2)(1+uw^2)+tu^2w(1-u^2)^{-2}(A^{(1)}-uA^{(3)})(A^{(3)}-uA^{(1)})\right)\sum_{q=0}^\infty v_{0,q,0}u^q}\\
  & \displaystyle{-tu(1+uw^2)\sum_{q=0}^\infty v_{0,q,1}u^q-wu(1+ut^2)\sum_{q=0}^\infty v_{1,q,0}u^q-(tw+tu^2w(1-u^2)^{-1})v_{0,0,0}}.
\end{array}
\end{equation}
From Equations (\ref{R1}) and (\ref{R3}), we also deduce
\begin{equation}\label{R7b}
    \sum_{q=0}^\infty v_{1,q,0}u^q=(1-u^2)^{-1}(A^{(1)}-uA^{(3)})\sum_{q=0}^\infty v_{0,q,0}u^q.
\end{equation}
So, by using Equations (\ref{R7}) and (\ref{R7b}), we obtain
\begin{equation}\label{rel4bis}
\begin{array}{rl}
 & \displaystyle{(1-tA^{(1)}+t^2A^{(2)}-t^3A^{(3)}+t^4)(1-wA^{(3)}+w^2A^{(2)}-w^3A^{(1)}+w^4)x}\\ \\
 = & \displaystyle{\Big((1+ut^2)(1+uw^2)-tu(1+uw^2)(1-u^2)^{-1}(A^{(3)}-uA^{(1)})}\\
  & \displaystyle{-wu(1+ut^2)(1-u^2)^{-1}(A^{(1)}-uA^{(3)})+tu^2w(1-u^2)^{-2}(A^{(1)}-uA^{(3)})(A^{(3)}-uA^{(1)})\Big)\sum_{q=0}^\infty v_{0,q,0}u^q}\\
  & \displaystyle{-(tw+tu^2w(1-u^2)^{-1})v_{0,0,0}},
\end{array}
\end{equation}
i.e., by multiplying (\ref{rel4bis}) by $(1-u^2)^2$,
\begin{equation}\label{rel4ter}
\begin{array}{rl}
 & \displaystyle{(1-u^2)^2(1-tA^{(1)}+t^2A^{(2)}-t^3A^{(3)}+t^4)(1-wA^{(3)}+w^2A^{(2)}-w^3A^{(1)}+w^4)x}\\ \\
 = & \displaystyle{\Big((1-u^2)^2(1+ut^2)(1+uw^2)-tu(1+uw^2)(1-u^2)(A^{(3)}-uA^{(1)})}\\
  & \displaystyle{-wu(1+ut^2)(1-u^2)(A^{(1)}-uA^{(3)})+tu^2w(A^{(1)}-uA^{(3)})(A^{(3)}-uA^{(1)})\Big)\sum_{q=0}^\infty v_{0,q,0}u^q}\\
  & \displaystyle{-(tw(1-u^2)^2+tu^2w(1-u^2))v_{0,0,0}}.
\end{array}
\end{equation}
$\bullet$ Consider now the following equation
\begin{equation}\label{R2}
    A^{(2)}\sum_{q=0}^\infty v_{0,q,0}u^q=u^{-1}\sum_{q=0}^\infty v_{0,q,0}u^q+u\sum_{q=0}^\infty v_{0,q,0}u^q+u\sum_{q=0}^\infty v_{1,q,1}u^q-u^{-1}v_{0,0,0}.
\end{equation}
Then, according to Equation (\ref{R8}), we have
$$\begin{array}{rcl}
\displaystyle{A^{(2)}\sum_{q=0}^\infty v_{0,q,0}u^q} & = & \displaystyle{u^{-1}\sum_{q=0}^\infty v_{0,q,0}u^q+u\sum_{q=0}^\infty v_{0,q,0}u^q}\\
     & & \displaystyle{+u(1-u^2)^{-2}(A^{(1)}-uA^{(3)})(A^{(3)}-uA^{(1)})\sum_{q=0}^\infty v_{0,q,0}u^q-u(1-u^2)^{-1}v_{0,0,0}-u^{-1}v_{0,0,0}},
  \end{array}
$$
i.e.
\begin{equation}\label{R9}
    \left(A^{(2)}-u^{-1}-u-u(1-u^2)^{-2}(A^{(1)}-uA^{(3)})(A^{(3)}-uA^{(1)})\right)\sum_{q=0}^\infty v_{0,q,0}u^q=-(u(1-u^2)^{-1}+u^{-1})v_{0,0,0}.
\end{equation}
This last equation reads
\begin{equation}\label{R9bis}
    \left(-u(1-u^2)^2A^{(2)}+(1+u^2)(1-u^2)^2+u^2(A^{(1)}-uA^{(3)})(A^{(3)}-uA^{(1)})\right)\sum_{q=0}^\infty v_{0,q,0}u^q=(1-u^2)v_{0,0,0}.
\end{equation}
Now, by using Equations (\ref{rel4ter}) and (\ref{R9bis}), we get
\begin{equation}\label{rel5}
\begin{array}{rl}
 & \displaystyle{(1-u^2)^2(1-tA^{(1)}+t^2A^{(2)}-t^3A^{(3)}+t^4)(1-wA^{(3)}+w^2A^{(2)}-w^3A^{(1)}+w^4)}\\
  & \displaystyle{\left(-u(1-u^2)^2A^{(2)}+(1+u^2)(1-u^2)^2+u^2(A^{(1)}-uA^{(3)})(A^{(3)}-uA^{(1)})\right)x}\\ \\
 =  & \displaystyle{-tw(1-u^2)\left(-u(1-u^2)^2A^{(2)}+(1+u^2)(1-u^2)^2+u^2(A^{(1)}-uA^{(3)})(A^{(3)}-uA^{(1)})\right)v_{0,0,0}}\\
 & \displaystyle{\Big((1-u^2)^2(1+ut^2)(1+uw^2)-tu(1+uw^2)(1-u^2)(A^{(3)}-uA^{(1)})}\\
  & \displaystyle{-wu(1+ut^2)(1-u^2)(A^{(1)}-uA^{(3)})+tu^2w(A^{(1)}-uA^{(3)})(A^{(3)}-uA^{(1)})\Big)(1-u^2) v_{0,0,0}},\\
\end{array}
\end{equation}
i.e., after simplification by $(1-u^2)$,
\begin{equation}\label{rel5bis}
\begin{array}{rl}
 & \displaystyle{(1-tA^{(1)}+t^2A^{(2)}-t^3A^{(3)}+t^4)(1-wA^{(3)}+w^2A^{(2)}-w^3A^{(1)}+w^4)}\\
  & \displaystyle{\left((1+u^2)(1-u^2)^2-u(1-u^2)^2A^{(2)}+u^2(A^{(1)}-uA^{(3)})(A^{(3)}-uA^{(1)})\right)x}\\ \\
 =  & \displaystyle{\Big((1-u^2)((1+ut^2)(1+uw^2)-tw(1+u^2))+twu(1-u^2)A^{(2)}}\\
  & \displaystyle{-tu(1+uw^2)(A^{(3)}-uA^{(1)})-wu(1+ut^2)(A^{(1)}-uA^{(3)})\Big) v_{0,0,0}}.\\
\end{array}
\end{equation}
The proposition is proved. $\blacksquare$

\subsection{\textsf{An inversion formula}}

\noindent In order to inverse the relation obtained in Proposition \ref{form} and get an explicit expression for $P_\Gamma(t,\,u)$, we need the rational function $f$ defined by
$$\begin{array}{rlcl}
  f\ : & \mathbb{C}^3 & \rightarrow & \mathbb{C}(t,\,u,\,w) \\
   & (d_1,\,d_2,\,d_3) & \mapsto & \displaystyle{(1-td_1+t^2d_2-t^3d_3+t^4)^{-1}(1-wd_3+w^2d_2-w^3d_1+w^4)^{-1}}\\
 & & & \displaystyle{\Big((1+u^2)(1-u^2)^2-u(1-u^2)^2d_2+u^2(d_1-ud_3)(d_3-ud_1)\Big)^{-1}}. \\
\end{array}$$
According to Proposition \ref{form}, we may write
 $$\begin{array}{rcl}
J(t,\,u,\,w)\,v_{0,0,0} & = &  T_\Gamma\left(1-t\Lambda^{(1)}+t^2\Lambda^{(2)}-t^3\Lambda^{(3)}+t^4\right)\left(1-w\Lambda^{(3)}+w^2\Lambda^{(2)}-w^3\Lambda^{(1)}+w^4\right) \\
 & & \!\!\!\!\!\!\!\!\!\!\!\!\Big((1+u^2)(1-u^2)^2-u(1-u^^2)^2\Lambda^{(2)}+u^2(\Lambda^{(1)}-u\Lambda^{(3)})(\Lambda^{(3)}-u\Lambda^{(1)})\Big)T_\Gamma^{-1}P_\Gamma(t,\,u,\,w).
   \end{array}$$
We deduce that
\begin{equation}\label{fond}
P_\Gamma(t,\,u,\,w)=T_\Gamma\,\Delta(t,\,u,\,w)\,T_\Gamma^{-1}J(t,\,u,\,w)\,v_{0,0,0}=(T_\Gamma\,\Delta(t,\,u,\,w)\,T_\Gamma)\,(T_\Gamma^{-2}J(t,\,u,\,w)\,v_{0,0,0}),
\end{equation}
where $\Delta(t,\,u,\,w)\in\mathbf{M}_{l+1}\mathbb{C}(t,\,u,\,w)$ is the diagonal matrix defined by $$\Delta(t,\,u,\,w)_{jj}=f(\Lambda^{(1)}_{jj},\,\Lambda^{(2)}_{jj},\,\Lambda^{(3)}_{jj})=f\left(\overline{\chi(g_j)},\,\frac{1}{2}(\chi(g_j)^2-\chi(g_j^2)),\,\chi(g_j)\right).$$

\noindent This last formula proves Theorem \ref{mainth}.\\

\begin{Rk}$\\$
The Poincar\'e series $\widehat{P}_\Gamma(t)$ of the algebra of invariants $\mathbb{C}[z_1,\,z_2,\,z_3,\,z_4]^\Gamma$ is given by $$\widehat{P}_\Gamma(t)=P_\Gamma(t,\,0,\,0)_0=P_\Gamma(0,\,0,\,t)_0.$$
\end{Rk}

\subsection{\textsf{Remark for $\mathbf{SL}_n\mathbb{C}$}}

\noindent In this section, we consider an integer $n\geq 2$ and a subgroup $\Gamma$ of $\mathbf{SL}_n\mathbb{C}$. As in paragraph \ref{firstsection}, let $\{\gamma_0,\dots,\,\gamma_l\}$ be the set of equivalence classes of irreducible finite dimensional complex representations of $\Gamma$, where $\gamma_0$ is the trivial
representation. The character associated to $\gamma_j$ is denoted by $\chi_j$.\\
Consider $\gamma : \Gamma\rightarrow\mathbf{SL}_n\mathbb{C}$ the natural $n-$dimensional representation, and $\chi$ its character.
By complete reducibility we get the decomposition $\gamma_j\otimes\gamma=\bigoplus_{i=0}^l a^{(1)}_{ij}\gamma_i$ for every $j\in\llbr 0,\,l\rrbr$, and we set
$A^{(1)}:=\left(a^{(1)}_{ij}\right)_{(i,j)\in\llbr 0,\, l\rrbr^2}\in\mathbf{M}_{l+1}\mathbb{N}$.\\

\noindent Let $\go{h}$ be a Cartan subalgebra of $\go{sl}_n\mathbb{C}$ and let $\varpi_1,\dots,\ \varpi_{n-1}$ be the corresponding fundamental weights, and $V(p_1\varpi_1+\dots+p_{n-1}\varpi_{n-1})$ the simple $\go{sl}_n\mathbb{C}-$module of highest weight $p_1\varpi_1+\dots+p_{n-1}\varpi_{n-1}$ with $\mathbf{p}:=(p_1,\dots\,p_{n-1})\in \mathbb{N}^{n-1}$. Then we get an irreducible representation
$\pi_{\mathbf{p}}:\mathbf{SL}_n\mathbb{C}\rightarrow\mathbf{GL}(V(p_1\varpi_1+\dots+p_{n-1}\varpi_{n-1}))$.
The restriction of $\pi_{\mathbf{p}}$ to the subgroup $\Gamma$ is a representation of~$\Gamma$, and by complete reducibility, we get the decomposition ${\pi_{\mathbf{p}}}|_{\Gamma}=\bigoplus_{i=0}^l m_i(\mathbf{p})\gamma_i$, where the $m_i(\mathbf{p})$'s are non negative integers.
Let $\mathcal{E}:=(e_0,\dots,\,e_l)$ be the canonical basis of
$\mathbb{C}^{l+1}$, and
$$\displaystyle{v_{\mathbf{p}}:=\sum_{i=0}^l m_i(\mathbf{p})e_i\in \mathbb{C}^{l+1}}.$$\\
As $\gamma_0$ is the trivial representation, we have $v_{\mathbf{0}}=e_0$. Let us consider the vector (with elements of $\mathbb{C}\llbr t_1,\dots,\,t_{n-1}\rrbr=\mathbb{C}\llbr \mathbf{t}\rrbr$ as coefficients)
$$\displaystyle{P_\Gamma(\mathbf{t}):=\sum_{\mathbf{p}\in \mathbb{N}^{n-1}}v_{\mathbf{p}}\mathbf{t}^\mathbf{p}\in\left(\mathbb{C}\llbr \mathbf{t}\rrbr\right)^{l+1}},$$
and denote by $\displaystyle{P_\Gamma(\mathbf{t})_j}$ its $j-$th coordinate in the basis $\mathcal{E}$.\\

\noindent Given the results from Kostant (\cite{Kos85} and \cite{Kos06}) for $\mathbf{SL}_2\mathbb{C}$ and our results (\cite{BP09}) about $\mathbf{SL}_3\mathbb{C}$, we then formulate the following conjecture:

\begin{Conj}$\\$
The coefficients of the vector $P_\Gamma(\mathbf{t})$ are rational fractions in $\mathbf{t}$, i.e. the formal power series $P_\Gamma(\mathbf{t})_i$ are rational functions $$P_\Gamma(\mathbf{t})_i:=\frac{N_\Gamma(\mathbf{t})_i}{D_\Gamma(\mathbf{t})},\ i\in\llbr 0,\,l\rrbr,$$
where the $N_\Gamma(\mathbf{t})_i$'s and $D_\Gamma(\mathbf{t})$ are elements of $\mathbb{Q}[\mathbf{t}]$.
\end{Conj}

\section{\textsf{An example of explicit computation}}

\noindent The classification of finite subgroups of $\mathbf{SL}_4\mathbb{C}$ is given in \cite{HH01}. It consists in infinite series and $30$ exceptional groups (types $I,\,II,\dots,\,XXX$). We give here an explicit computation of $P_\Gamma(t,\,u,\,w)$ for one of these exceptional groups.
Consider the matrices
$$F_1=\left( \begin {array}{cccc} 1&0&0&0\\\noalign{\medskip}0&1&0&0\\\noalign{\medskip}0&0&j&0\\\noalign{\medskip}0&0&0&{j}^{2}
\end {array} \right),\ F_2'=\frac{1}{3}\left( \begin {array}{cccc} 3&0&0&0\\\noalign{\medskip}0&-1&2&2\\\noalign{\medskip}0&2&-1&2\\\noalign{\medskip}0&2&2&-1\end {array}
 \right),\ F_3'=\frac{1}{4}\left( \begin {array}{cccc} -1&{\sqrt{15}}&0&0\\\noalign{\medskip}{\sqrt{15}}&1&0&0\\\noalign{\medskip}0&0&0&4
\\\noalign{\medskip}0&0&4&0\end {array} \right),$$
and the subgroup $\Gamma=\langle F_1,\ F_2',\ F_3'\rangle$ of $\mathbf{SL}_4\mathbb{C}$ (type $II$ in \cite{HH01}).\\
Here $l=4$,
$$A^{(1)}=A^{(3)}=\left( \begin {array}{ccccc} 0&0&0&1&0\\\noalign{\medskip}0&0&1&1&1\\\noalign{\medskip}0&1&0&1&1\\\noalign{\medskip}1&1&1&1&1
\\\noalign{\medskip}0&1&1&1&2\end {array} \right)\ \ \textrm{and}\ \ A^{(2)}= \left( \begin {array}{ccccc} 0&1&1&0&0\\\noalign{\medskip}1&1&0&1&2
\\\noalign{\medskip}1&0&1&1&2\\\noalign{\medskip}0&1&1&2&2
\\\noalign{\medskip}0&2&2&2&2\end {array} \right).
$$\\
$\textrm{rank}(A^{(1)})=\textrm{rank}(A^{(2)})=4$, and the eigenvalues of $A^{(1)},\ A^{(2)},\ A^{(3)}$ are\\
\begin{small}
$\Theta^{(1)}=\overline{\Theta^{(3)}}=( 4,\ 0,\ -1,\ 1,\ -1 ),\ \Theta^{(2)}=(6,\ -2,\ 1,\ 0,\,1)$, $p=4$, and $\tau_0=s_0s_1$, $\tau_1=s_2$, $\tau_2=s_3$, $\tau_3=s_4$.\\
\end{small}

\noindent According to formula \ref{fond}, we get
$$\begin{array}{rcl}
D_{\Gamma}(t,\,u,\,w) & = & \left( w-1 \right) ^{4} \left( u+1 \right) ^{3} \left( u-1 \right) ^{5} \left( t-1 \right) ^{4} \left( {w}^{2}+w+1 \right)  \left( {w}^{4}+
{w}^{3}+{w}^{2}+w+1 \right)\\
 & & \left( w+1 \right) ^{2} \left( {u}^{4}+{u
}^{3}+{u}^{2}+u+1 \right)  \left( {u}^{2}+u+1 \right) ^{2} \left( {t}^
{2}+t+1 \right)  \left( {t}^{4}+{t}^{3}+{t}^{2}+t+1 \right)  \left( t+
1 \right) ^{2}\\
 & = & (u-1)(u+1)(u^2+u+1)\widetilde{D}_{\Gamma}(t)\widetilde{D}_{\Gamma}(u)\widetilde{D}_{\Gamma}(w),
\end{array}$$\\
with $\widetilde{D}_{\Gamma}(t)=(t-1)^4(t+1)^2(t^2+t+1)(t^4+t^3+t^2+t+1)$. Moreover,
$$\widehat{P}_{\Gamma}(t)={\frac {{t}^{8}-{t}^{6}+{t}^{4}-{t}^{2}+1}{{t}^{12}-2\,{t}^{10}-{t}^{9
}+{t}^{8}+{t}^{7}+{t}^{5}+{t}^{4}-{t}^{3}-2\,{t}^{2}+1}}.$$\\

\noindent Because of the to big size of the numerators $N_{\Gamma}(t,\,u,\,w)_i$'s, only the denominator is given in the text: all the numerators may be found on the web (\verb"http://math.univ-lyon1.fr/~butin/").

\begin{small}

\end{small}

\end{document}